\documentclass[a4paper,11pt]{amsart}
\usepackage[all]{xy}
\usepackage{amscd}
\usepackage{amsfonts}
\usepackage{a4wide}
\usepackage{amsmath, amsthm}
\usepackage{amssymb}
\usepackage{fullpage}
\usepackage[applemac]{inputenc}
\usepackage[pdftex,breaklinks,colorlinks,citecolor=blue]{hyperref}

\newtheorem{theorem}{Theorem}

\newtheorem{lemma}[theorem]{Lemma}
\newtheorem{proposition}[theorem]{Proposition}
\newtheorem{remark}[theorem]{Remark}

\title{Enumeration of a special class of irreducible polynomials in characteristic 2}

\author{Alp Bassa  $^\P$}
\address{Alp Bassa. Bo\u{g}azi\c{c}i University,
Faculty of Arts and Sciences,
Department of Mathematics,
34342 Bebek, \.{I}stanbul,
Turkey,
}
\email{alp.bassa@boun.edu.tr}
\author{Ricardo Menares $^\dag$}
\address{Ricardo Menares. Pontificia Universidad Cat\'olica de Chile,  Facultad de Matem\'aticas, Vicu\~na Mackenna 4860, Santiago, Chile.}
\email{rmenares@mat.uc.cl}

\thanks{$\P$ Bo\u gazi\c ci University}
\thanks{\dag Pontificia Universidad Cat\'olica de Chile}
\thanks{Both authors were partially supported by Conicyt-MEC 80130064 grant. Alp Bassa was partially supported by the BAGEP Award of the Science Academy with funding supplied by Mehve{\c s} Demiren in memory of Selim Demiren and by Bo\u gazi\c ci University Research Fund Grant Number 15B06SUP2. Ricardo Menares was partially supported by  FONDECYT 1171329  grant. }

\begin{document}

\begin{abstract}
$A$-polynomials were introduced by Meyn and play an important role in the iterative construction of high degree self-reciprocal irreducible polynomials over the field $\mathbb F_2$, since they constitute the starting point of the iteration. The exact number of $A$-polynomials of each degree was given by Niederreiter. Kyuregyan extended the construction of Meyn to arbitrary even finite fields. We relate the $A$-polynomials in this more general setting to inert places in a certain extension of elliptic function fields and obtain an explicit counting formula for their number. In particular, we are able to show that, except for an isolated exception, there exist A-polynomials of every degree. 
\end{abstract} 
\maketitle

The $Q$-transform plays a prominent role in the construction of (self-reciprocal) irreducible polynomials. Given a polynomial $f\in \mathbb F_q[T]$, its $Q$-transform is given by
$$f^Q(T):=T^{\deg f}\cdot f\Bigl(T+\frac1T\Bigr).$$
Then, $f^Q$ is a self-reciprocal polynomial of degree $2\cdot\deg f$. Clearly for $f^Q$ to be irreducible a necessary condition is that $f$ is irreducible. In characteristic 2, a simple sufficient and necessary condition for the irreducibility of $f^Q$ in terms of the coefficients of $f$ was established by Meyn (\cite{meyn}, Theorem 6). More surprising is the fact that it is even possible to devise criteria to ensure that irreducibility is preserved under repeated applications of the $Q$-transform, hence giving an infinite sequence of self-reciprocal irreducible polynomials of increasing degree. Starting with an irreducible polynomial $f_0\in \mathbb F_q[T]$, we   iteratively define a sequence of polynomials $f_m \in \mathbb F_q[T]$ by
$$f_{m+1}=f_m^Q,  \quad m \geq 0. $$

In what follows, we set $q=2^r$ and describe the conditions mentioned above in characteristic 2. Let $f(T)=T^n+a_{n-1}T^{n-1}+\ldots+a_1T+a_0 \in \mathbb F_q[T]$ be a monic irreducible polynomial of degree $n$. We say that $f$ is an $A$-\emph{polynomial} if $Tr_{\mathbb F_q/\mathbb F_2}(a_{n-1})=1$ and $Tr_{\mathbb F_q/\mathbb F_2}(a_1/a_0)=1$.  Then, whenever $f_0$ is an $A$-polynomial, the polynomial $f_m(T)$ is irreducible of degree $n2^m$ for all $m$.  This fact was proved by  Meyn \cite{meyn} and Varshamov \cite{var}, when $r=1$, and generalized latter by Kyuregyan \cite{recurrent} for general $r \geq 1$. 

When $q$ is odd, the $Q$-transform behaves in a more subtle way and the results are less complete. In that setting, S. Cohen  introduced the related $R$-transform, which leads to a comparable iterative construction of irreducible reciprocal polynomials \cite{SCohen}. For a comparison of both transforms in terms of Galois theory see \cite{GaloisTheory}. 

In this note, we provide a closed formula for the number of $A$-polynomials in $\mathbb F_{2^r}[T]$ of given degree. On the other hand, it is known that there are no $A$-polynomials of degree 3 in $\mathbb F_2[T]$. In all other cases, our formula allows us to deduce the existence of $A$-polynomials of each degree over every finite field of characteristic 2.

\newpage

\begin{theorem}\label{main}
We denote by $A_r(n)$  the number of $A$-polynomials in $\mathbb F_{2^r}[T]$ of degree $n$.

\begin{enumerate}
 \item[a)] Write $n=2^k\cdot m$, ($m$ odd) and let $\alpha=\frac{-1+\sqrt{-7}}{2}$. Then, 
$$A_r(n)=\frac{1}{4n} \sum_{d|m}\mu\Bigl(\frac{ m}{d}\Bigr)\bigl(q^{2^kd}+1-\alpha^{r 2^kd}-{\overline{\alpha}}^{r 2^k d}\bigr).$$
Here, $\mu$ is the M\"obius function. In particular, we have that 
\begin{equation}\label{estimate}
\left|A_r(n)-\frac{q^n}{4n}\right| \leq \frac{\sigma_0(m)}{4n}(q^{n/3}+1+2\cdot 2^{rn/6}).
\end{equation}
 Here, $\sigma_0(m)$ is the number of positive divisors of $m$. 
  Hence,  $A_r(n) \sim \frac{q^n}{4n}$ as $ n \rightarrow \infty.$

\item[b)] Assume $(r,n)\neq (1,3)$. Then, $A_r(n)\geq 1$. 

\end{enumerate}
\end{theorem}

In order to prove Theorem \ref{main},  we exploit the correspondence between irreducible polynomials of each degree $n$ and  degree $n$ places that are inert in a particular unramified extension of elliptic function fields. Then, we  show that an exact count can be obtained using the corresponding $L$-polynomials. A very rough estimate then ensures the existence  of $A$-polynomials of any degree $n\geq 7$ over all even finite fields, which can be used as the starting point of an iterative construction. 
  
Theorem \ref{main} is a generalization of a result of Niederreiter in \cite{niederreiter}, treating the case $r=1$ (cf. Remark \ref{nied}). His method requires an explicit evaluation of certain Kloosterman sums attached to additive characters, which is available only when the base field is small. By turning around Niederreiter's reasoning, we obtain an archimedean evaluation of a certain weighted average of Kloosterman sums (cf. Proposition \ref{Kloos} below).

For results and notation about algebraic function fields we refer the reader to \cite{sti}. In \cite{MoisioRanto} a related approach is applied to the problem of counting polynomials with prescribed coefficients.

\section{Interpreting $A$-polynomials in terms of  an extension of elliptic function fields}
Let $q=2^r$, and let $\mathbb F_q$  be the finite field of $q$ elements. Consider the rational function field $F=\mathbb F_q(x)$ and the extensions $E_1=F(y_1)$ and $E_2=F(y_2)$, with $y_1^2+y_1=x$ and $y_2^2+y_2=1/x$. Both $E_1$ and $E_2$ are again rational function fields.

\begin{lemma}
Let $f(x)=x^n+a_{n-1}x^{n-1}+\ldots+a_1x+a_0 \in \mathbb F_q[x]$ be a monic irreducible polynomial of degree $n$, with $f(x)\neq x$.  We denote by  $P_f$  the place of $F$ of degree $n$ associated to $f$. Then,

\begin{enumerate}
 \item[i)] $P_f$ is inert in $E_1/F$ if and only if  $Tr_{\mathbb F_q/\mathbb F_2}(a_{n-1})=1$
 \item[ii)] $P_f$ is inert in $E_2/F$ if and only if  $Tr_{\mathbb F_q/\mathbb F_2}(a_1/a_0)=1.$

\end{enumerate}
 In particular, $f$ is an $A$-polynomial if and only if $P_f$ is inert in both extensions $E_1/F$ and $E_2/F$.
\end{lemma}

\proof First we prove i). If $c$ is a root of $f$ in $\overline{\mathbb F}_q$, then $c \neq 0$ and $a_{n-1}=Tr_{\mathbb F_{q^n}/\mathbb F_q}(c)$. Hence by the transitivity of the trace, the condition $Tr_{\mathbb F_q/\mathbb F_2}(a_{n-1})=1$ is equivalent to $Tr_{\mathbb F_{q^n}/\mathbb F_2}(c)\neq 0$. By Hilbert's Theorem 90, this happens exactly if $c$ is not of the form $\gamma^2-\gamma$ for any $\gamma\in \mathbb F_{q^n}$. In turn, this happens if and only if $f(y_1^2+y_1) \in \mathbb F[y_1]$ is irreducible. The last condition is equivalent to  $P_f$ being inert in $E_1/F$, thus proving i).

Part ii) follows form Part i) applied to $f^*(x)=x^nf(1/x)\quad  \diamond$\\

Consider the compositum $E'=E_1\cdot E_2$ over $F$.  The extension $E'/F$ is Galois, with Galois group $\mathbb Z/2\mathbb Z\times \mathbb Z/2\mathbb Z$, with $E_1$ and $E_2$ corresponding to the subgroups $\mathbb Z/2\mathbb Z\times \{0\}$ and $\{0\}\times \mathbb Z/2\mathbb Z$, respectively. Let $E$ be the subfield corresponding to third subgroup $H$, the diagonal subgroup. Clearly $E=F(y)$ with $y=y_1+y_2$ satisfying  $y^2+y=x+1/x$. In other words, $E$ is the function field of the elliptic curve over $\mathbb F_q$ with $j=1$. 

\begin{center}
\mbox{\xymatrix{
&E'=E_1\cdot E_2\ar@{-}[dl] \ar@{-}[dr] \ar@{-}^{H}[d]&\\
E_1=F(y_1) \ar@{-}[dr] & E=F(y_1+y_2) \ar@{-}[d] & E_2=F(y_2) \ar@{-}[dl] \\
& F=\mathbb F_q(x)&\\
}}
\end{center}

\begin{proposition}\label{connection}
Let $P'$ be a place of $E'$ above $P_f$. We denote by $G(P'|P_f)$ be the associated decomposition group in $E'/F$. Let $C_r(n)$ be the number if inert places of degree $n$ in the extension $E'/E$. Then,

\begin{enumerate}
 \item[i)] $f$ is an $A$-polynomial if and only if $G(P'|P_f)=H$
 \item[ii)] We have that $C_r(n)=2A_r(n)$.  
\end{enumerate}
\end{proposition}

\proof 
The polynomial $T$ corresponds to the zero $P_T$ of $x$. Since only    the pole $P_\infty$ of $x$   ramifies  in the extension $E_1/F$, the places $P_T$ and $P_\infty$ are the only places of $F$ ramified in the extension $E'/F$. Both places are ramified in $E/F$ and in each case the place of $E$ lying above them splits in the extension $E'/E$. The extension $E'/E$ is unramified. Using the Riemann-Hurwitz genus formula, we see that the genera of $E$ and $E'$ are both $1$. For a place $P_f\neq P_T, P_\infty$, let $Z(P_f)$ be the associated decomposition group in the extension $E'/F$. The place $P_f$ is inert in $E_1/F$ and $E_2/F$, if and only if $Z(P_f)=H$.  Hence to each degree $n$ place $P_f$ of $F$ that is inert in $E_1/F$ and $E_2/F$ there correspond two places of $E$ of degree $n$ that are inert in $E'/E$. 
In particular, to every $A$-polynomial $f$  there correspond two places of $E$ of degree $\deg f$ that are inert in $E'/E$.  $\quad  \diamond$

\section{Enumerating $A$-polynomials over arbitrary even finite fields}
In this section we provide a proof of Theorem \ref{main}. Following Proposition \ref{connection}, we need to count the number of inert places of $E$ and $E'$.  We will achieve this task by means of their $L$-polynomials.

 All function fields can be defined already over $\mathbb F_2$. Hence we consider the extension $\mathbb F_2(x,y)/\mathbb F_2(x)$ with $y^2+y=x+1/x$. Among the rational places of $\mathbb F_2(x)$, the pole and zero of $x$ are ramified, the zero of $x-1$ splits in $\mathbb F_2(x,y)/\mathbb F_2(x)$, giving a total of $4$ rational places of $\mathbb F_2(x,y)$. The elliptic function field hence has trace $-1$ and $L$-polynomial 
$$2t^2+t+1=(1-\alpha t)(1-{\overline{\alpha}} t) \text{ with } \alpha=\frac{-1+\sqrt{-7}}{2}.$$
The $L$-polynomial of the constant field extension $E=\mathbb F_q(x,y)$ with $q=2^r$ is hence 
$$L_E(t)=(1-\alpha^r t)(1-{\overline{\alpha}}^r t).$$

The $L$-polynomial $L_{E'}$ of $E'$ has to be divisible by $L_E$ and be also of degree $2$, since $g(E')=1$. Hence $L_{E'}=L_{E}$. The number of degree $n$ places for each of the function fields is given by (see \cite[Propositions 5.1.16 and 5.2.9]{sti})
\begin{equation}
B(n)=\frac{1}{n}\sum_{d|n}\mu\Bigl(\frac{n}{d}\Bigr) \bigl(q^d+1-\alpha^{r d}-{\overline{\alpha}}^{r d}\bigr).\label{Bn}
\end{equation}
For even $n$ the $B(n)$ places of $E'$ of degree $n$ come from the $C_r({n/2})$ inert places of degree $n/2$ of $E$ and the $B(n)-C_r(n)$ splitting places of degree $n$ (we get $2$ degree $n$ places for each splitting place). For odd $n$ the $B(n)$ places of $E'$ of degree $n$ come only from the $B(n)-C_r(n)$ splitting places of $E$ of degree $n$. Hence we obtain
\begin{eqnarray*}
B(n)=&C_r({n/2})+2\cdot (B(n)-C_r(n))& \text{for}\ n\ \text{even}\\
B(n)=&2\cdot (B(n)-C_r(n))& \text{for}\ n\ \text{odd}.\\
\end{eqnarray*}
Writing $n=2^k\cdot m$ for an odd integer $m$, we obtain
\begin{equation}\label{Cn}
C_r(2^k\cdot m)=\sum_{i=1}^{k+1}\frac{1}{2^i}B({2^{k+1-i}\cdot m}).
\end{equation}

Using Proposition \ref{connection}, ii) and equations~\eqref{Bn} and \eqref{Cn}, we obtain the formula stated in the first part of Theorem \ref{main}. The estimate \eqref{estimate} is obtained by using that $|\alpha|=\sqrt{2}$ and the fact that any proper divisor $d |m$ satisfies $d \leq m/3$. 

The function field $E$ has genus $1$. Hence using estimates for the number of higher degree places (see for instance\cite[Corollary 5.2.10]{sti}) $B(n)\geq 1$ for any $n$ with  $q^{(n-1)/2}(q^{1/2}-1)\geq 3$. Hence for $r\geq 4$ we have $B(n)\geq 1$ for any $n$. For $r=1$, we need $n\geq 7$ and for $r=2$ or $3$ we need $n\geq 3$ to ensure $B(n)\geq 1$ using this estimate. A case by case analysis using the $L$-polynomial shows that in all cases except $(r,n)=(1,3),$ we have $B(n)\geq 1$. 
Since $E'$ has genus 1, using  Proposition \ref{connection}, ii)  and equation~\eqref{Cn}, we obtain the second assertion in Theorem \ref{main}.

\begin{remark}\label{nied}
 Note that for $r=1$ we recover the main Theorem in \cite{niederreiter}. Indeed, use the elementary identity $$\alpha^t+\overline{\alpha}^t=\frac{1}{2^{t-1}}\sum_{j=0}^{[t/2]}\left(\begin{array}{c}
t\\
2j
\end{array}\right)(-1)^{t+j}7^j,$$ valid for all integers $t\geq 1$.
\end{remark}

\section{Averages of Kloosterman sums}

Let $q=2^r$ and let  $\chi : \mathbb F_q  \rightarrow \mathbb C^*$ be an additive character. For any integer $n\geq 1$, let $\chi^{(n)}  :\mathbb F_{q^n}  \rightarrow \mathbb C^*$ be the additive character defined by 
$$\chi^{(n)}(u)=\chi\circ Tr_{\mathbb F_{q^n}/\mathbb F_{2}}(u).$$

Let $$K( \chi^{(n)};a,b)=\sum_{\alpha \in \mathbb F_{q^n}^*} \chi^{(n)}(a\alpha + b\alpha^{-1})$$ 
be the Kloosterman sum attached to $\chi^{(n)}$.

\begin{proposition}\label{Kloos}
Assume $\chi$ is a non trivial additive character. Then,
\begin{equation}\label{archimedean}
\frac{1}{q^2}\sum_{u \in \mathbb F_q}\chi(u)\sum_{a,b \in \mathbb F_q \atop a+b=u} K(\chi^{(n)};a,b) = \sum_{d|n}\mu \left( \frac{n}{d}\right) \cdot d \cdot A_r(d),
\end{equation} for all $n \geq 1$.
\end{proposition}

\proof let 

\begin{align*}
R(n) =&\{ \alpha \in \mathbb F_{q^n}^* : Tr_{\mathbb F_{q^n}/\mathbb F_{2}}(\alpha)=Tr_{\mathbb F_{q^n}/\mathbb F_{2}}(\alpha^{-1})=1 \}.\\
R^*(n)=&\{\alpha\in R(n) : [\mathbb F_q(\alpha):\mathbb F_q]=n\}.
\end{align*}

We have that 
\begin{equation}\label{obvio}
|R^*(n)|=n\cdot A_r(n), \quad |R(n)| = \sum_{d|n} |R^*(d)|.
\end{equation}

Since $\chi$ is non trivial, for all $u \in \mathbb F_q$ we have that \cite[Corollary 5.31]{LN},  $$\frac{1}{q} \sum_{a \in \mathbb F_q}\chi(ua) = \left\{ \begin{array}{cc}
1 & \textrm{ if } u=0\\
0 & \textrm{ otherwise.}
\end{array}\right.$$ 

 Set $T^{(n)}:=Tr_{\mathbb F_{q^n}/\mathbb F_{2}}.$ We have that

\begin{align*}
|R(n)| &= \sum_{\alpha \in \mathbb F_{q^n}^*} 1_{\{T^{(n)}(\alpha)=1\}}\cdot  1_{\{T^{(n)}(\alpha^{-1})=1\}}\\
&=  \sum_{\alpha \in \mathbb F_{q^n}^*} \frac{1}{q^2} \sum_{a \in \mathbb F_q}\chi\left((T^{(n)}(\alpha)+1)a\right) \sum_{b \in \mathbb F_q}\chi\left((T^{(n)}(\alpha^{-1})+1)b\right)\\
&= \frac{1}{q^2} \sum_{\alpha \in \mathbb F_{q^n}^*} \sum_{a,b \in \mathbb F_q}\chi(a+b)\chi\left(T^{(n)}(a\alpha+b\alpha^{-1})\right)\\
&= \frac{1}{q^2} \sum_{u\in \mathbb F_q} \chi(u) \sum_{a,b \in \mathbb F_q \atop a+b=u} \sum_{\alpha \in \mathbb F_{q^n}^*}\chi^{(n)}(a\alpha+b\alpha^{-1})\\
&= \frac{1}{q^2}\sum_{u \in \mathbb F_q}\chi(u)\sum_{a,b \in \mathbb F_q \atop a+b=u} K(\chi^{(n)};a,b).
\end{align*}

We conclude by combining this equality with the relations \eqref{obvio}. $\quad  \diamond$

\begin{remark}
When $r=1$, Niederreiter finds in \cite{niederreiter} an explicit evaluation in archimedean terms of each individual Kloosterman sum in the LHS of \eqref{archimedean}. It seems difficult to obtain a similar individual evaluation for general $r$.  Our method proceeds by interpreting the whole LHS in terms of the number of places of appropriate degree in a particular elliptic extension and then use the explicit knowledge of the corresponding $L$-polynomials to determine such quantities. 
\end{remark}

\end{document}